\documentclass[a4paper,leqno,12pt,twoside]{amsart}
\usepackage{amsmath,amscd,amssymb,upgreek} 
\usepackage[all]{xy}
\usepackage{graphicx}
\usepackage{ifthen}
\usepackage{yhmath}
\usepackage{pdfsync}
\usepackage{color}
\usepackage{hyperref}

\usepackage[applemac] {inputenc}

\newcommand{\application}[5]{
\ifthenelse{\equal{#1}{0}}{}{#1:}\begin{array}[t]{ccl}
	#2 & \longrightarrow & #3
	\ifthenelse{\equal{#4}{0}}{}{ \\ #4 & \longmapsto & #5} 
\end{array}}
\newcommand{\Exemple}[1]{\begin{Exple}\emph{#1}\end{Exple}}
\newcommand{\Remarque}[1]{\begin{Rem}\emph{#1}\end{Rem}}
\newcommand{\Theoreme}[1]{\begin{The}#1\end{The}}
\newcommand{\Proposition}[1]{\begin{Prop}#1\end{Prop}}
\newcommand{\Lemme}[1]{\begin{Lem}#1\end{Lem}}

\newcommand{\Definition}[1]{\begin{Def}#1\end{Def}}
\newcommand{\Preuve}[1]{\noindent \textbf{Proof:} #1~\hfill$\blacksquare$\\}

\def\ord{\operatorname{ord}}
\def\R{\mathbb{R}}
\def\C{\mathbb{C}}

\def\N{\mathbb{N}}

\def\Z{\mathbb{Z}}

\def\S{\mathbb{S}}
\def\P{\mathbb{P}}

\def\Aff{\operatorname{Aff}}

\def\Mc{\operatorname{\mathcal{M}c}}
\def\SEP{\operatorname{SEP}}
\def\LVMB{\operatorname{LVMB}}
\def\LVM{\operatorname{LVM}}
\def\Pn{\P^{n}(\C)}

\begin{document}
\newtheorem{The}{Theorem}[section]
\newtheorem{Prop}[The]{Proposition}
\newtheorem{Def}[The]{Definition}
\newtheorem{Lem}[The]{Lemma}
\newtheorem{Cor}[The]{Corollary}
\newtheorem{Rem}[The]{Remark}
\newtheorem{Exple}[The]{Example}

\begin{title}
{LVMB manifolds and quotients of toric varieties}
\end{title}
\begin{author}
{L. Battisti}
\end{author}
\begin{abstract}
In this article, we study a class of manifolds introduced by Bosio (see~\cite{Bosio:2001aa}) called $\LVMB$ manifolds. We provide an interpretation of his construction in terms of quotient of toric manifolds by complex Lie groups. Furthermore, $\LVMB$ manifolds extend a class of manifolds obtained by Meersseman in \cite{Meersseman:2000aa}, called $\LVM$ manifolds, and we give a characterization of these manifolds using our toric description. Finally, we give an answer to a question asked by Cupit-Foutou and Zaffran in \cite{Cupit-Foutou:2007aa}.
\end{abstract}
\thanks{\noindent \!\!\!\!\!\!\!Laurent BATTISTI - Aix-Marseille Universit\'e - Laboratoire d'Analyse, Topologie, Probabilités UMR~7353~CNRS, 39 Rue Joliot-Curie, 13013 Marseille, France.\\Financial support for this PhD thesis is assured by the French \emph{Région Provence-Alpes-Côte d'Azur}. E-mail: \texttt{battisti@math.cnrs.fr}}

\maketitle 

\tableofcontents

\setcounter{section}{-1}
\section{Introduction}
In \cite{Lopez-de-Medrano:1997aa}, López de Medrano and Verjovsky introduce a family of complex compact manifolds, obtained as quotients of a dense open subset $U$ of $\Pn$ by the action of a complex Lie group isomorphic to $\C$. This construction is  extended to the case of an action of $\C^{m}$ (with $m$ a positive integer) by Meersseman in \cite{Meersseman:2000aa}, and these manifolds are called $\LVM$ manifolds.

Then, Bosio extends in \cite{Bosio:2001aa} the construction due to Meersseman by allowing other actions of $\C^{m}$ on certain open subsets of $\Pn$, and these manifolds are called $\LVMB$ manifolds. In short, given a family $\mathcal{E}_{m,n}$ of subsets of $\{0, ..., n\}$ of cardinal $2m+1$ (we assume $n$ and $m$ are integers such that $2m\leqslant n$) and a family $\mathcal{L}$ of $n+1$ linear forms on $\C^{m}$ satisfying some conditions, Bosio associates to $\mathcal{E}_{m,n}$ an open subset $U$ of $\Pn$ and to $\mathcal{L}$ an action of $\C^{m}$ on $\Pn$, such that the quotient $U/\C^{m}$ is a compact complex manifold. We say that the pair $(\mathcal{E}_{m,n},\mathcal{L})$ is an $\LVMB$ datum, and that it is an $\LVM$ datum if the manifold we obtain is an $\LVM$ manifold. We will recall this construction in details later. ~\\

Our first goal in this article is to express Bosio's construction in terms of toric geometry. We find a fan $\Delta$ in $\R^{n}$ such that the corresponding toric manifold is the open set $U$ and we see how the projection of this fan by a suitable $2m$-dimensional linear subspace of $\R^{n}$ helps to understand the action of $\C^{m}$ on $U$. Roughly speaking, the set $\mathcal{E}_{m,n}$ will correspond to the fan $\Delta$, and the choice of the linear forms on $\C^{m}$ will give the subspace $\R^{2m}$. The converse will also work, i.e. with suitable conditions on a fan $\Delta$ and the choice of a suitable $2m$-dimensional subspace of $\R^{n}$, we will get an $\LVMB$ datum. We have:~\\

\noindent\textbf{Theorem~\ref{BosioIsToric}.}\emph{
i) Let $(\mathcal{L},\mathcal{E}_{n,m})$ be an $\LVMB$ datum. Then there is a pair $(E,\Delta)$ where $E$ is a $2m$-dimensional linear subspace of $\R^{n}$ and $\Delta$ is a subfan of the fan $\Delta_{\Pn}$ of $\Pn$, satisfying the following two properties:
\begin{itemize}
\item[$a)$] the projection map $\pi: \R^{n} \rightarrow \R^{n}/E\cong\R^{n-2m}$ is injective on $|\Delta|$,
\item[$b)$] the fan $\pi(\Delta)$ is complete in $\R^{n}/E$, i.e. $|\pi(\Delta)| = \R^{n}/E$.
\end{itemize}
ii) Conversely, given a pair $(E,\Delta)$ having the two properties above, one obtains an $\LVMB$ datum.}~\\

To prove this theorem, we first need to slightly extend the notion of the \emph{manifold with corners} associated to a fan. This is what we do in the first section.~\\

In \cite{Bosio:2001aa}, Bosio gives a criterion for deciding whether an $\LVMB$ datum is an $\LVM$ datum or not. Our second goal is to translate this criterion in our new toric setting, it is given by the following:~\\

\noindent\textbf{Theorem~\ref{lvmbEquivPolytopal}.}\emph{
Let $(\mathcal{L},\mathcal{E}_{m,n})$ be an $\LVMB$ datum and $(E,\Delta)$ its associated pair given by theorem~\ref{BosioIsToric}. Then, $(\mathcal{L},\mathcal{E}_{m,n})$ is an $\LVM$ datum if and only if the projection by $E$ of the fan $\Delta$ is polytopal.}~\\

Finally, we will use this characterization to show that if two $\LVMB$ manifolds $X$ and $Y$ are biholomorphic and if $X$ is an $\LVM$ manifold, then $Y$ is also an $\LVM$ manifold. This question was outlined and partly answered by Cupit-Foutou and Zaffran in \cite{Cupit-Foutou:2007aa} and we answer it with theorem~\ref{questionByZaffran}:~\\

\noindent\textbf{Theorem~\ref{questionByZaffran}.}\emph{
Let $(\mathcal{L}_{1}, \mathcal{E}_{m,n})$ and $(\mathcal{L}_{2}, \mathcal{E}'_{m',n'})$ be two $\LVMB$ data giving two biholomorphic $\LVMB$ manifolds, then $n=n'$, $m=m'$ and $(\mathcal{L}_{1}, \mathcal{E}_{m,n})$ is an $\LVM$ datum if and only if $(\mathcal{L}_{2}, \mathcal{E}'_{m,n})$ is an $\LVM$ datum.}~\\

This article is organized as follows: in the first part, we extend the notion of the manifold with corners of a fan and we give some basic properties of it. Here, the fan we consider is non-necessarily rational. We then use this new object in the second part, where we prove theorem~\ref{BosioIsToric} after recalling Bosio's construction. We also study the case when two $\LVMB$ data give the same pair $(E,\Delta)$. This leads us to a correspondence statement between the set of all such pairs and the quotient of the set of $\LVMB$ data by a suitable equivalence relation. In the third part we detect $\LVM$ data among $\LVMB$ data by proving theorem~\ref{lvmbEquivPolytopal}, and in the fourth part we use this criterion to obtain theorem~\ref{questionByZaffran}.

\section{The manifold with corners of a fan}
The main definition of this section (the \emph{manifold with corners} of a fan and its topology) is a basic one in the theory of toric manifolds, and it can be found in \cite{Ash:2010aa} and \cite{Oda:1988aa} for instance. Here we extend this definition to a fan without the ``rationality condition'', that is, the vectors generating the cones we work with are not necessarily located in some rational lattice of a vector space. We also give basic properties of the topology of the associated manifold with corners that will be of use later. ~\\

\subsection{Definition of the manifold with corners associated to a fan}
Let $N_{\R}$ be a real vector space of dimension $n$. 
\Definition{A subset $\sigma$ of $N_{\R}$ is called a \textbf{convex polyhedral cone (with apex at the origin)} if there exist a finite number of elements $v_{1}, ..., v_{m}$ such that 
\[\sigma = \R_{\geqslant0}v_{1}+ ... + \R_{\geqslant 0} v_{m},\] and that $\sigma \cap (-\sigma) = \{0\}$, where $\R_{\geqslant0}$ is the set of non-negative real numbers. The \textbf{dimension} of a cone is the dimension of the smallest linear subspace of $N_{\R}$ containing this cone. We denote by $L(\sigma)$ this vector space. We say that $\sigma$ is \textbf{simplicial} if the vectors generating this cone are linearly independent.}

\Definition{The cone in $N_{\R}^{*}$ (the dual of $N_{\R}$) \textbf{dual} to $\sigma$ is the set 
$$\check{\sigma}:=\{\varphi \in N_{\R}^{*}~|~\varphi(x)\geqslant 0~\text{ for all }x\in \sigma\}.$$  It is a convex polyhedral cone. A subset $\tau$ of $\sigma$ is a \textbf{face} of $\sigma$ if there is a $\varphi_{0}\in\check{\sigma}$ such that $$\tau=\{x \in \sigma~|~\varphi_{0}(x) = 0\}.$$ We denote this relation between $\tau$ and $\sigma$ by $\tau < \sigma$. The cone $\sigma$ is a face of itself, and we say that $\tau < \sigma$ is a \textbf{proper face} of $\sigma$ if $\tau \neq \sigma$.}

\Remarque{In the following we always, for short, say ``\textbf{cone}'' for a convex polyhedral cone with apex at the origin, since we will only consider such sets.}

\Definition{A \textbf{fan} $\Delta$ of $N_{\R}$ is a set of cones satisfying the following two properties: \begin{itemize}
\item[-] each face of a cone of $\Delta$ is a cone of $\Delta$,
\item[-] the intersection of two cones of $\Delta$ is a face of each of these cones.
\end{itemize}
The \textbf{support} of a fan $\Delta$ is the set
$$|\Delta|:= \bigcup_{\sigma\in\Delta}\sigma.$$}


Let $\Delta$ be a fan of $N_{\R}$. For each cone $\sigma\in\Delta$, $L(\sigma)$ denotes the linear subspace of $N_{\R}$ generated by $\sigma$ and we denote by $N_{\R}^{\sigma}$ any complementary subspace of $L(\sigma)$ in $N_{\R}$ (that is, $N_{\R} =  N_{\R}^{\sigma} \oplus L(\sigma)$). Let $(x_{n})_{n\in\N}$ be a sequence of points in $N_{\R}$ and write the decomposition $x_{n}=y_{n}+z_{n} \in N_{\R}^{\sigma}\oplus L(\sigma)$ for each $n \in \N$.

\Definition{\label{defNDelta}Given a fan $\Delta$, we call $\boldsymbol{\mathcal{N}_{\Delta}}$ the set of sequences $(x_{n})_{n\in\N}$ of points in $N_{\R}$ such that there exists a cone $\sigma((x_{n})_{n\in\N}) \in \Delta$ satisfying: 
\begin{itemize}
\item[-] there is a point $y \in N_{\R}^{\sigma}$ such that $\displaystyle \lim_{n \rightarrow \infty }y_{n} = y$ and
\item[-] for every $w \in L(\sigma)$, there is an integer $p$ such that $z_{n}\in w + \sigma$ for all $n\geqslant p$.
\end{itemize}}

There are two remarks to be made here:
\Remarque{The set $\mathcal{N}_{\Delta}$ is well-defined, that is, its definition does not depend on the choice of a complementary subspace $N_{\R}^{\sigma}$ of $L(\sigma)$ for each $\sigma \in \Delta$. Given a cone $\sigma$, let $N_{\R}^{\sigma}$ and $N_{\R}^{'\sigma}$ be two complementary subspaces of $L(\sigma)$, and write the two decompositions $x_{n}=y_{n}+z_{n}\in N_{\R}^{\sigma}\oplus L(\sigma)$ and $x_{n}=y'_{n}+z'_{n}\in N_{\R}^{'\sigma}\oplus L(\sigma)$ for a sequence $(x_{n})_{n\in\N}$ of $N_{\R}$. It is straightforward to check that if the sequences $(y_{n})_{n\in \N}$ and $(z_{n})_{n\in \N}$ satisfy the conditions of the previous definition, then $(y'_{n})_{n\in \N}$ and $(z'_{n})_{n\in \N}$ also do.}

\Remarque{\label{rqBornee}If $(x_{n})_{n\in\N} \in \mathcal{N}_{\Delta}$ is a bounded sequence of $N_{\R}$, then necessarily $\sigma((x_{n})_{n\in\N})=\{0\}$.}

The following lemma shows the uniqueness of $\sigma \in \Delta$ satisfying the conditions of the previous definition.

\Lemme{\label{uniciteConeMc}Let $(x_{n})_{n\in\N}\in \mathcal{N}_{\Delta}$. Then, the cone $\sigma((x_{n})_{n\in\N})$ is unique. }
\Preuve{
Assume that there exist two distinct cones $\sigma, \sigma' \in \Delta$ both satisfying the conditions of definition~\ref{defNDelta}. Their intersection $\tau$ is a face of each of these cones and we can assume that it is different from at least one of them, say $\tau\neq\sigma$, otherwise $\tau=\sigma=\sigma'$.  

Write the two decompositions $x_{n} = y_{n}+z_{n}=y_{n}'+z_{n}'$ for all $n\in \N$ corresponding to $\sigma$ and $\sigma'$ respectively. The condition on the sequence $(x_{n})_{n\in\N}$ for $\sigma$ means that $(y_{n})_{n\in\N}$ converges, therefore it is bounded; thus there is a compact set $K_{\sigma}$ and an integer $p$ such that the sequence $(x_{n})_{n\in\N}$ has values in $K_{\sigma} + \sigma$ for $n\geqslant p$; similarly there is a compact set $K_{\sigma'}$ and an integer $p'$ such that the elements of the sequence $(x_{n})_{n\in\N}$ all lie in $K_{\sigma'} + \sigma'$ for $n\geqslant p'$. 

As a consequence, there is a compact set $K$ such that all the elements of $(x_{n})_{n\in\N}$ are (for $n\geqslant \max(p,p')$) in the intersection of $K+\sigma$ and $K+\sigma'$. 
Chose $w$ in a complementary subspace $V$ of $L(\tau)$ in $L(\sigma)$ such that $w+\sigma$ and $\sigma'$ have empty intersection. It is possible because $\tau$ is a proper face of $\sigma$, so $\dim L(\tau) +1 \leqslant  \dim L(\sigma)$. Now choose $\lambda>0$ large enough such that $K+\lambda w + \sigma$ and $K+\sigma'$ also have empty intersection. By definition of $\mathcal{N}_{\Delta}$, there is an integer $p''$ such that the elements of $(x_{n})_{n\in\N}$ are also all in $K + \lambda w + \sigma$ and $K + \sigma'$ for $n\geqslant p''$. This is a contradiction.}

Now, we define an equivalence relation on $\mathcal{N}_{\Delta}$. Keeping the previous notations, one has the following definition:
\Definition{Two sequences $(x_{n})_{n\in\N}$ and $(x'_{n})_{n\in\N}$ in $\mathcal{N}_{\Delta}$ are equivalent when $\sigma((x_n)_{n\in\N}) = \sigma((x'_n)_{n\in\N})$ and $y-y'\in L(\sigma)$. In this case, we shall write $(x_{n})_{n\in\N}\sim(x'_{n})_{n\in\N}$ and denote by $y + \infty \cdot \sigma$ the equivalence class of $(x_{n})_{n\in\N}$. \\
The \textbf{manifold with corners of $\boldsymbol{\Delta}$} is the quotient of $\mathcal{N}_{\Delta}$ by this equivalence relation. We denote this space by $\boldsymbol{\operatorname{\mathcal{M}c}(\Delta)}$.}

To every $x\in N_{\R}$ and each cone $\sigma\in\Delta$, we associate an element of $\operatorname{\mathcal{M}c}(\Delta)$, denoted by $x+\infty\cdot\sigma$. It is the equivalence class of any sequence $(x_{n})_{n\in\N}\in\mathcal{N}_{\Delta}$ such that $y=p(x)$ where $p(x)$ is the projection of $x$ on a complementary subspace of $L(\sigma)$. 

We use the conventions $x+\infty\cdot\{0\}=x$ and $x+\infty\cdot \sigma = x' + \infty\cdot \sigma$ for two points $x,x' \in N_{\R}$ with $x-x'\in L(\sigma)$.

\Remarque{The vector space $N_\R$ is a subset of the manifold with corners of $\Delta$. For $x \in N_\R$, we simply consider the constant sequence $(x_n)_{n\in \N}$ with $x_{n}=x$. We still denote by $N_{\R}$ the image of $N_{\R}$ in $\Mc(\Delta) =\mathcal{N}_{\Delta}/\sim$ and specify $N_{\R}\subset \Mc(\Delta)$ if necessary. Notice that any cone $\sigma$, as a subset of $N_{\R}$, is then also a subset of $\Mc(\Delta)$.}

\Remarque{As one would expect, if $\Delta$ is a rational fan and $X_{\Delta}$ is the associated toric variety, the set $\Mc(\Delta)$ is homeomorphic to the quotient $X_{\Delta} / (\S^{1})^{n}$, i.e. the manifold with corners of $X_{\Delta}$. For this fact, we refer to \cite{Ash:2010aa}, section~I.1.}
~\\
\subsection{The topology of a manifold with corners}~\\
Now that the manifold with corners of a fan is defined, we endow this space with a topology and we then prove some properties of this topology that we will use later. These properties are very close to what we have in the case of toric manifolds: there is an action of $N_{\R}$ on the manifold with corners of a fan, and it is compact if and only if the fan is finite and complete.~\\

As before, let $\Delta$ be a fan in a vector space $N_{\R}$.
\Definition{\label{topoMc}We equip $\Mc(\Delta)$ with a topology in the following way: a neighbourhood basis of a point $y + \infty \cdot \sigma \in \Mc(\Delta)$ is the collection of sets \[U_{\varepsilon, w}(y + \infty \cdot \sigma):= \bigcup_{ \tau < \sigma}(y+w+B_{\varepsilon}+\sigma + \infty \cdot \tau),\] 
for $\varepsilon> 0$ and $w \in L(\sigma)$ ($B_{\varepsilon}$ is the unit open ball around $0$ in $N_{\R}$).}

Figure~\ref{fig1} below depicts a fan $\Delta=\{\{0\},\tau,\tau_{1}, \tau_{2}, \sigma\}$ in $\R^{2}$, its associated manifold with corners and the neighbourhood of a point.

\begin{figure}[h]
\includegraphics[scale=1]{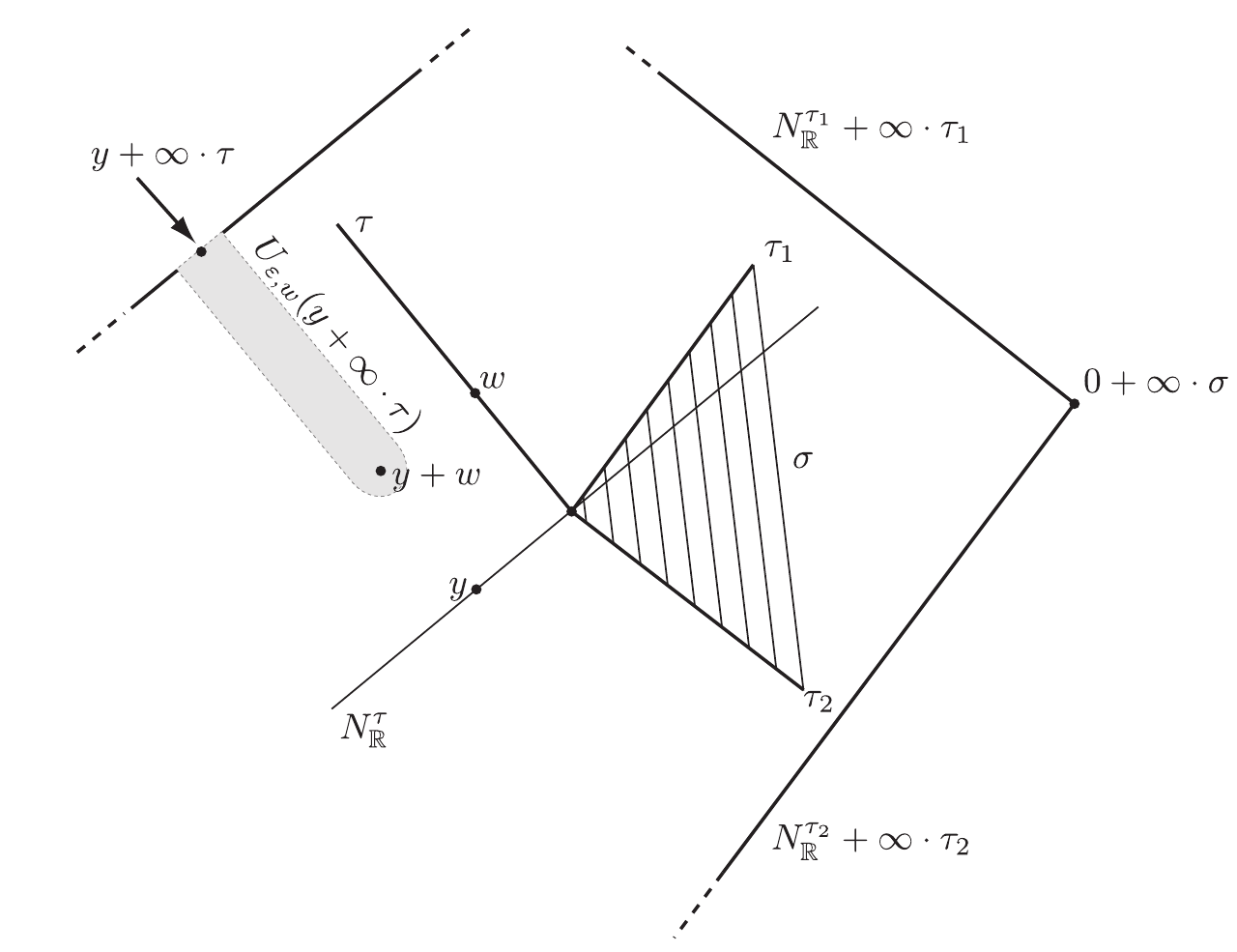}
\caption{\label{fig1}A fan and its manifold with corners}
\end{figure}

\Lemme{The space $\Mc(\Delta)$ has a natural continuous action of $N_{\R}$ which extends the action of $N_{\R}$ on itself.}
\Preuve{
Let $y + \infty\cdot \sigma \in \Mc(\Delta)$ (this point being $y$ if $\sigma = \{0\}$ or $0+\infty \cdot \sigma$ if $\sigma$ is $n$-dimensional), and $x \in N_{\R}$. Set $$x.(y + \infty\cdot\sigma):= p(x)+y + \infty \cdot \sigma,$$ where $p(x)$ is the projection of $x$ to the linear subspace $N_{\R}^{\sigma}$. This action is clearly continuous. }

\Remarque{\label{rqLimite}Let $(x_{n})_{n\in\N}$ be an element of $\mathcal{N}_{\Delta}$ and $y+\infty\cdot \sigma$ its equivalence class in $\Mc(\Delta)$. It is easy to see that if we consider $(x_{n})_{n\in\N}$ as a sequence of points in $\Mc(\Delta)$, then it converges to $y+\infty\cdot \sigma$. }

\Lemme{\label{coneCompact}Let $\sigma$ be a cone of a fan $\Delta$. Then the closure $S$ of $\sigma \subset \Mc(\Delta)$ in $ \Mc(\Delta)$ is compact.}
\Preuve{
We consider a sequence of $S$, denoted by $(x_{n})_{n\in\N}$. Call $\mathcal{S}:=S\setminus \sigma$ and write $S$ as the disjoint union $S= \sigma \sqcup \mathcal{S}$. We distinguish cases. \\

$\alpha$) First, assume that $(x_{n})_{n\in\N}$ has infinitely many elements in $\sigma$. After a possible extraction of subsequence, we can assume that the sequence $(x_{n})_{n\in\N}$ has all its elements in this set. Then, either $(x_{n})_{n\in\N}$ is bounded and we are done, or it is not bounded. In this case, after another extraction of subsequence if needed, we can assume that $(\|x_{n}\|)_{n\in\N}$ is strictly increasing, where $\|\cdot\|$ is a norm on $N_{\R}$. Denote by $\tau_{1}, ..., \tau_{r}$ the vectors generating $\sigma$, and write for all $n \in \N$: 
\begin{equation}\label{eqDecompKappa} x_{n} = x_{1,n}\tau_{1} + ... + x_{r,n}\tau_{r}.\end{equation}
Note that this decomposition is not necessarily unique (in fact it is only the case if $\sigma$ is simplicial). Since $(x_{n})_{n\in\N}$ is not bounded, one at least of the sequences $(x_{i,n})_{n\in\N}$ (for $i \in \{1, ..., r\}$) is also not bounded. Assume (to simplify writing) that the sequences $(x_{1,n})_{n\in\N}, ..., (x_{j,n})_{n\in\N}$ are bounded while the sequences $(x_{j+1,n})_{n\in\N}, ..., (x_{r,n})_{n\in\N}$ are not. Consider the cone generated by $\tau_{j+1}, ..., \tau_{r}$. We call this cone $\kappa$ if it is a proper face of $\sigma$, otherwise we set $\kappa:=\sigma$. The sequence $(x_{1,n})_{n\in\N}$ is bounded so it admits a convergent subsequence, say $(x_{1, \varphi(n)})_{n\in\N}$. Similarly, the sequence $(x_{2,\varphi(n)})_{n\in\N}$ is bounded, so we can find a convergent subsequence for it, and so on. After extracting subsequences a finite number of times, we can therefore assume that $(x_{1,n})_{n\in\N}, ..., (x_{j,n})_{n\in\N}$ are convergent and that $(x_{j+1,n})_{n\in\N}, ..., (x_{r,n})_{n\in\N}$ are strictly increasing. 

Set $y_{n}:=\pi(x_{1,n}\tau_{1}+...+ x_{j,n} \tau_{j})$ for all $n\in\N$ where $\pi$ is the projection on a complementary space $N_{\R}^{\kappa}$ of $L(\kappa)$ and denote by $y$ the limit of this sequence. Notice that if $\kappa=\sigma$, we have $y_{n}=0$ for all $n\in \N$. Write, for all $n\in \N$, $x_{n}=y_{n}+z_{n}$ with $z_{n}=x_{n}-y_{n}$. One sees that for all $n\in\N$, $z_{n}\in L(\kappa)$ hence the sequence $(x_{n})_{n\in\N}$ is an element of $\mathcal{N}_{\Delta}$. Indeed, the sequence $(y_{n})_{n\in\N}$ is convergent, taking values in a complementary space $N_{\R}^{\kappa}$ of $L(\kappa)$ and the sequence $(z_{n})_{n\in \N}$ satisfies, for all $w \in L(\kappa)$, the existence of a rank $p$ such that for all $n\geqslant p$, $z_{n} \in w+\kappa$. Lemma~\ref{uniciteConeMc} tells us that $\kappa$ is uniquely defined, hence the fact that the decomposition~($\ref{eqDecompKappa}$) above is not unique has no incidence. According to remark~\ref{rqLimite}, we have proved that the sequence $(x_{n})_{n\in\N}$ admits a convergent subsequence with limit $y+ \infty\cdot \kappa \in S$. \\

$\beta$) Now, assume that $(x_{n})_{n\in\N}$ has infinitely many terms in $\mathcal{S}$. As above, after a possible subsequence extraction we may assume that all the values of $(x_{n})_{n\in\N}$ are in $\mathcal{S}$. By the previous part of the proof, we see that every element of $\mathcal{S}$ is of the form $\pi(y)+ \infty\cdot \kappa$ where $y\in \sigma$, $\kappa$ is a face of $\sigma$ and $\pi$ is the projection on a complementary space $N_{\R}^{\kappa}$ of $L(\kappa)$ (in particular, $y=0$ if $\kappa=\sigma$). 
Since $\sigma$ only possesses a finite number of faces, after a subsequence extraction if needed, we may assume that there is a unique face $\kappa$ of $\sigma$ such that for all $n \in \N$, $x_{n}=y_{n}+\infty\cdot \kappa$ with $y_{n} \in \pi(\sigma)$. If $\kappa=\sigma$, the sequence is constant hence the result is proved. Now, assume that $\kappa$ is a proper face of $\sigma$. As earlier, we denote by $\tau_{1}, ..., \tau_{r}$ the generators of $\sigma$, the face $\kappa$ being generated over $\R_{\geqslant0}$ by $\tau_{j+1}, ..., \tau_{r}$. 
We now use the same reasoning as in the first part of the proof, this time applied to the sequence \[y_{n} = y_{1,n}\pi(\tau_{1}) + ... + y_{j,n}\pi(\tau_{j})\]
(notice that since $\kappa$ is a face of $\sigma$, $\pi(\sigma)$ is a cone in $N_{\R}^{\kappa}$). As before, there is an integer $\ell \in\{1,..., j+1\}$ such that the sequences $(y_{1,n})_{n\in\N}, ..., (y_{\ell,n})_{n\in\N}$ are convergent (none of them is convergent if $\ell=j+1$ by convention), $(y_{\ell+1,n})_{n\in\N}, ..., (y_{j,n})_{n\in\N}$ are strictly increasing and we see that $(x_{n})_{n\in\N}$ converges in $\Mc(\Delta)$ to  a point $y'+\infty\cdot \kappa'$ where $\kappa'$ is either the proper face of $\sigma$ generated by $\tau_{\ell}, ..., \tau_{r}$, either $\sigma$ itself (and $\kappa$ is a face of $\kappa'$ in each case). 
%
This concludes the proof.
}

We now prove the following proposition, which just extends what we already now in the toric case:
\Proposition{\label{compactComplet}The manifold with corners associated to a fan $\Delta$ in $\R^{n}$ is compact if and only if $\Delta$ is finite and complete.}
\Preuve{If we assume that $\Delta$ is finite and complete, the result is a consequence of lemma~\ref{coneCompact}. \\
Assume now that $\Mc(\Delta)$ is compact. First, suppose $\Delta$ is not finite and take any sequence $(x_{n})_{n\in\N}$ of points in $\Mc(\Delta)$ with $x_{n} \in N_{\R}+\infty\cdot \sigma_{n}$ and $\sigma_{n}\neq\sigma_{p}$ for $n\neq p$. The definition of the topology of $\Mc(\Delta)$ tells us that if a subsequence of $(x_{n})_{n\in\N}$ were convergent, then a subsequence of $(\sigma_{n})_{n\in\N}$ would become constant after some rank, which is a contradiction.

Now that we know $\Delta$ is finite, assume it is not complete. There is a vector $v\in N_{\R}$ such that $\R_{\geqslant0} v \subset \R^{n}\setminus|\Delta|$. We now claim that the sequence $(nv)_{n\in\N}$ has no convergent subsequence in $\Mc(\Delta)$. If it were the case, its limit would lie in $\Mc(\Delta)\setminus N_{\R}$, i.e. there would exist $\sigma \in \Delta$ and $y\in N_{\R}^{\sigma}$ such that $\displaystyle \lim_{n\rightarrow +\infty}nv = y+\infty\cdot \sigma$. The definition of the topology tells us that for $\varepsilon>0$, there exists an integer $p$ such that for $n\geqslant p$, $nv \in \displaystyle \bigcup_{\tau < \sigma}(y+B_{\varepsilon}+\sigma + \infty\cdot \tau)$ , i.e. $nv \in y+B_{\varepsilon}+\sigma$. This is impossible because for all $y\in N_{\R}^{\sigma}$ and $\varepsilon>0$, the set $(y+B_{\varepsilon}+\sigma)\cap \R_{\geqslant 0}v$ is bounded. Indeed, it is a convex set (it is an intersection of convex sets) therefore connected, which means that it is an interval of $\R_{\geqslant0}v$. If it were not bounded, this would mean that $\R_{\geqslant a}v\subset y+B_{\varepsilon}+\sigma$ for some $a\geqslant 0$. Then, one would have $av \in (y+B_{\varepsilon}+\sigma) \setminus \sigma$ (because $v\not\in\sigma$) and we obtain a contradiction because every point of $(y+B_{\varepsilon}+\sigma)\setminus\sigma$ is mapped outside $y+B_{\varepsilon}+\sigma$ by homotheties of center $0$ and large enough ratio. 
We thus have proved that $(nv)_{n\in\N}$ has no convergent subsequence and $\Mc(\Delta)$ is not compact if $\Delta$ is not complete.}

\section{Describing $\LVMB$ manifolds in terms of toric geometry}
In this section we translate Bosio's construction in terms of toric geometry. In Bosio's construction, one needs an open subset $U$ of $\Pn$ and an $m$-dimensional complex Lie group $G$ (with $2m\leqslant n$) such that $U/G$ is a compact complex manifold. First we recall this construction, then we see how we can associate a fan $\Delta$ to $U$, a linear subspace $E$ of $\R^{n}$ to $G$, and see how the manifold with corners of the fan $\Delta$ helps understanding the quotient map $U \rightarrow U/G$.
\subsection{\label{BosioConstr}Bosio's construction}
The construction we explain here is due to Bosio and it generalizes a work of Meersseman (see~\cite{Bosio:2001aa} and \cite{Meersseman:2000aa} respectively).\\

Let $m,n$ be positive integers such that $2m \leqslant n$; let $\mathcal{L}:=(\ell_{0}, ..., \ell_{n})$ be $n+1$ linear forms of $\C^{m}$ such that any subfamily of $2m+1$ elements of $\mathcal{L}$ is an $\R$-affine basis of $(\C^{m})^{*}$, where $(\C^{m})^{*}$ is the dual space of $\C^{m}$. Call $\mathcal{E}_{m,n}$ a family of subsets of $\{0,..., n\}$, each having $2m+1$ elements. For every $P\in\mathcal{E}_{m,n}$, call $\mathcal{L}_{P}$ the corresponding subfamily of $\mathcal{L}$. We are interested in the two following conditions on $\mathcal{L}$ and $\mathcal{E}_{m,n}$:
\begin{itemize}
\item if for all $P \in \mathcal{E}_{m,n}$ and for all $i\in \{0, ..., n\}$, there exists $j \in P$ such that $(P\setminus\{j\})\cup \{i\} \in \mathcal{E}_{m,n}$, we say that $(\mathcal{L}, \mathcal{E}_{m,n})$ satisfies the $\boldsymbol{\SEP}$ (for \textbf{s}ubstitute \textbf{e}xistence \textbf{p}rinciple),
\item if for all $P,Q \in \mathcal{E}_{m,n}$, the interiors of the convex envelopes of $\mathcal{L}_{P}$ and $\mathcal{L}_{Q}$ have non-empty intersection, we say that $(\mathcal{L}, \mathcal{E}_{m,n})$ satisfies the \textbf{imbrication condition}.
\end{itemize}

An $\boldsymbol{\LVMB}$ \textbf{datum} is a pair $(\mathcal{L}, \mathcal{E}_{m,n})$ satisfying these two conditions. Following Bosio's denomination, we say that an integer $i\in\{0, ..., n\}$ is \textbf{indispensable} if $i \in P$ for all $P\in \mathcal{E}_{m,n}$. If $i$ is indispensable, we also say that $\ell_{i}$ is indispensable.\\

Given a pair $(\mathcal{L},\mathcal{E}_{m,n})$ (not necessarily an $\LVMB$ datum), we construct an open subset of $\Pn$ and an action on this set by a Lie group $G$ and see when the quotient is a compact complex manifold:

First, call $U$ the (open) set of points $z=[z_{0}:...:z_{n}]\in \Pn$ such that there exists $P_{z}\in\mathcal{E}_{m,n}$ satisfying: for all $i\in P_{z}, z_{i}\neq 0$. 

Then, define an action of $G\cong\C^{m}$ on $U$ by:
$$\begin{array}{rcl}
\C^{m}\times U & \longrightarrow & U\\
(Z, [z_{0}:...:z_{n}])& \longmapsto & [\exp(\ell_{0}(Z))z_{0}:...:\exp(\ell_{n}(Z))z_{n}].\end{array}$$

Bosio relates the properness and the cocompactness of this action on $U$ with the $\SEP$ and the imbrication condition, this is the following
\Lemme{\label{lemmeBosio}The action of $G$ on $U$ is proper if and only if $(\mathcal{L},\mathcal{E}_{m,n})$ satisfies the imbrication condition.\\
If this action is proper, it is cocompact if and only if $(\mathcal{L},\mathcal{E}_{m,n})$ satisfies the $\SEP$.}

Hence the quotient $X_{n,m}:=U/G$ is a compact complex manifold if and only if $(\mathcal{L},\mathcal{E}_{m,n})$ is an $\LVMB$ datum. In this case, the manifold $X_{n,m}$ is of complex dimension $n-m$ and is called an $\boldsymbol{\LVMB}$ \textbf{manifold}.

~\\
\subsection{The toric viewpoint}~\\
As it is well-known, the complex projective space $\Pn$ is a toric manifold given by the fan $\Delta_{\Pn}$ in $\R^{n}$ which consists of the $n+1$ cones generated by $n$ of the $n+1$ vectors $e_{1}, ..., e_{n}, -(e_{1}+ ...+ e_{n})$ (where $(e_{1}, ..., e_{n})$ is the canonical basis of $\R^{n}$) and their faces. We say that a fan $\Delta$ is a \textbf{subfan} of $\Delta_{\Pn}$ if every cone $\sigma \in\Delta$ is a cone in $\Delta_{\Pn}$.\\

Our goal in this section is to prove the following
\Theoreme{\label{BosioIsToric}
i) Let $(\mathcal{L},\mathcal{E}_{n,m})$ be an $\LVMB$ datum. Then there is a pair $(E,\Delta)$ where $E$ is a $2m$-dimensional linear subspace of $\R^{n}$ and $\Delta$ is a subfan of the fan $\Delta_{\Pn}$ of $\Pn$, satisfying the following two properties:
\begin{itemize}
\item[$a)$] the projection map $\pi: \R^{n} \rightarrow \R^{n}/E\cong\R^{n-2m}$ is injective on $|\Delta|$,
\item[$b)$] the fan $\pi(\Delta)$ is complete in $\R^{n}/E$, i.e. $|\pi(\Delta)| = \R^{n}/E$.
\end{itemize}
ii) Conversely, given a pair $(E,\Delta)$ having the two properties above, one obtains an $\LVMB$ datum.}

\subsubsection{Preliminary results}
We prove the two parts of this theorem separately in the next two subsections. Before this, we need some preliminary statements: we see how the two conditions of the first part of this theorem are related with the $\SEP$ and the imbrication condition.~\\

The following two lemmas are respectively translations of the first and the second part of lemma~\ref{lemmeBosio} in our toric setting, in the sense that they characterize the properness and the cocompactness of an action of a linear space $E$ on a manifold with corners with conditions on its fan:
\Lemme{[Properness]\label{equivCond2}Let $\Delta$ be a fan in $\R^{n}$ and $E\cong \R^{k}$ be a linear subspace of $\R^{n}$. Then there is an action of $E$ on $\Mc(\Delta)$ and this action is proper if and only if the restriction of the quotient map $\pi: \R^{n} \rightarrow \R^{n}/E$ to the support $|\Delta|$ of $\Delta$ is injective.}
\Preuve{
First, assume that $\pi$ is not injective on $|\Delta|$, i.e. there exist two cones  $\sigma_{1}$ and $\sigma_{2}\in \Delta$, and $x \in \sigma_{1}$, $y\in \sigma_{2}$ such that $y-x \in E$. We know that the closures $K_{1}$ and $K_{2}$ in $\Mc(\Delta)$ of $\sigma_{1}$ and $\sigma_{2}$ respectively  are compact subsets (lemma~\ref{coneCompact}). We then write, for $\lambda \geqslant 1$, the equality $\lambda(y-x) + \lambda x = \lambda y$. Hence the set $\{t \in E ~|~ t K_{1} \cap K_{2}\neq \emptyset \}$ is not bounded, therefore the action of $E$ on $\Mc(\Delta)$ is not proper.~\\

Conversely, assume now that $E$ is not acting properly on $\Mc(\Delta)$. 
This means that there are two compact subsets $K_{1}$ and $K_{2}$ of $\Mc(\Delta)$ such that the set $\mathcal{E}:=\{t \in E ~|~ t K_{1} \cap K_{2}\neq\emptyset \}$ is not bounded. We introduce the following notation, for $\sigma \in \Delta$ and $i=1,2$: 
\[K_{i}^{\sigma}:=K_{i}\cap (\R^{n}+\infty \cdot \sigma).\]

Since the sets $K_{i}$ are compact, there exist only a finite number of cones $\sigma \in \Delta$ satisfying $K_{i}^{\sigma} \neq \emptyset$. We call $\Delta'$ the collection of these cones. 

We then have $K_{i} = \displaystyle\bigsqcup_{\sigma \in \Delta'} K_{i}^{\sigma}$. The action of $E$ induces an action on each component $\R^{n}+\infty\cdot \sigma$ for $\sigma \in \Delta$, so we have the following decomposition:
$$\mathcal{E} = \displaystyle\bigcup_{\sigma \in \Delta'}\{t \in E ~|~ t K_{1}^{\sigma} \cap K_{2}^{\sigma}\neq\emptyset \}.$$ Since the set $\Delta'$ is finite, there exists a cone $\sigma \in \Delta'$ such that $\mathcal{E}^{\sigma}:=\{t \in E ~|~ t K_{1}^{\sigma} \cap K_{2}^{\sigma}\neq\emptyset \}$ is not bounded. Consequently there is a sequence $(t_{n})_{n\in\N}$ of $\mathcal{E}^{\sigma}$ with $\|t_{n}\| \rightarrow + \infty$.

If $E \cap L(\sigma) \neq\{0\}$, we are done. Indeed, in this case there exists $x = \lambda_{1}\tau_{1}+ ...+ \lambda_{k}\tau_{k} \in E$ where the $\tau_{i}$'s are the generating rays of $\sigma$. 

We then write 
\begin{equation*}
x = \sum_{\lambda_{i}\geqslant 0}\lambda_{i}\tau_{i} - \sum_{\lambda_{j}<0}(-\lambda_{j})\tau_{j} = x^{+}-x^{-},
\end{equation*}

where the vectors $x^{+}$ and $x^{-}$ are both in $|\Delta|$, so the quotient map $\pi$ is not injective on $|\Delta|$.

Assume now that $E\cap L(\sigma) = \{0\}$. We now study $\R^{n}+\infty\cdot \sigma \cong \R^{n}/L(\sigma)=N_{\R}^{\sigma}$. 

Let $(x_{n})_{n\in\N}$ and $(y_{n})_{n\in\N}$ be sequences of $K_{1}^{\sigma}$ and $K_{2}^{\sigma}$ respectively (with at least one being not bounded) satisfying $t_{n}+x_{n}=y_{n}$. Up to extraction of subsequences, we may assume that they are convergent ($K_{1}$ and $K_{2}$ being compact sets), with respective limits $x_{0} + \infty\cdot \widetilde{\sigma_{1}}$ and $y_{0} + \infty\cdot \widetilde{\sigma_{2}}$ where $\sigma < \widetilde{\sigma_{i}}$ for $i=1,2$ (with eventually $\widetilde{\sigma_{i}}=\sigma$ for at most one index $i$). The definition of the topology of $\Mc(\Delta)$ then gives us that for any $\varepsilon>0$, there exists some $N$ such that for $n\geqslant N$, $x_{n} \in x_{0} + B_{\varepsilon} + \widetilde{\sigma_{1}}+\infty\cdot\sigma$ and $y_{n} \in y_{0} + B_{\varepsilon} + \widetilde{\sigma_{2}}+\infty\cdot\sigma$. 

We can then write $x_{n}=x_{0}+x_{\varepsilon,n}+x^{\widetilde{\sigma_{1}}}_{n}+\infty\cdot\sigma$ and $y_{n}=y_{0}+y_{\varepsilon,n}+y^{\widetilde{\sigma_{2}}}_{n}+\infty\cdot\sigma$ where $x_{\varepsilon,n}, y_{\varepsilon,n} \in B_{\varepsilon}$ and $x^{\widetilde{\sigma_{1}}}_{n}\in\widetilde{\sigma_{1}}$, $y^{\widetilde{\sigma_{2}}}_{n} \in \widetilde{\sigma_{2}}$.

The relation $t_{n}+x_{n}=y_{n}$ means
\begin{equation*}
t_{n}=p(t_{n}) = p(y_{0} - x_{0} + y_{\varepsilon,n} - x_{\varepsilon,n}) + p(y^{\widetilde{\sigma_{2}}}_{n} -  x^{\widetilde{\sigma_{1}}}_{n}),
\end{equation*}

where $p$ is the projection of $\R^{n}=E\oplus F\oplus L(\sigma)$ on $E\oplus F$ with respect to $L(\sigma)$.
Recall that $\|t_{n}\|\rightarrow + \infty$. We now write 
\begin{equation*}
\frac{t_{n}}{\|t_{n}\|} = \frac{p(y_{0} - x_{0} + y_{\varepsilon,n} - x_{\varepsilon,n})}{\|t_{n}\|} + \frac{p(y^{\widetilde{\sigma_{2}}}_{n} -  x^{\widetilde{\sigma_{1}}}_{n})}{\|t_{n}\|}.
\end{equation*}

First notice that the sequence $t_{n}/\|t_{n}\|$ is bounded, so we can assume it converges to some $t \in \R^{k}\setminus\{0\}$. Also notice that the sequence $(p(y_{0} - x_{0} + y_{\varepsilon,n} - x_{\varepsilon,n}))/\|t_{n}\|$ has limit $0$, therefore the sequence $(p(y^{\widetilde{\sigma_{2}}}_{n} -  x^{\widetilde{\sigma_{1}}}_{n}))/\|t_{n}\|$ also converges to $t$. Since the Minkowski sum of closed cones is closed, we know that $t\in p(\widetilde{\sigma_{2}}-\widetilde{\sigma_{1}})$. Therefore $t=y-x$ with $x\in p(\widetilde{\sigma_{1}})$ and $y\in p(\widetilde{\sigma_{2}})$. Moreover, there exists $z_{1}=\lambda_{1}\tau_{1}+...+\lambda_{k}\tau_{k} \in L(\sigma)$ such that $x+z_{1} \in \widetilde{\sigma_{1}}$. Then for all $z_{1}'=\mu_{1}\tau_{1}+...+\mu_{k}\tau_{k}$ with $\mu_{i} \geqslant \lambda_{i}$, we have $x+z_{1}' \in\widetilde{\sigma_{1}}$. Similarly, there exists $z_{2}\in L(\sigma)$ satisfying $y+z_{2}\in \widetilde{\sigma_{2}}$ and for all $z_{2}'=\mu'_{1} \tau_{1}+...+\mu'_{k}\tau_{k}$ with $\mu'_{i}\geqslant \lambda_{i}$, we have $y+z_{2}'\in\widetilde{\sigma_{2}}$. As a consequence there exists $w\in \sigma$ such that $y+w\in\widetilde{\sigma_{2}}$ and $x+w\in\widetilde{\sigma_{1}}$. This gives us $t=(y+w)-(x+w) \in \widetilde{\sigma_{2}}-\widetilde{\sigma_{1}}$, which implies the non-injectivity of $\pi$ on $|\Delta|$.}

For the second lemma, we need to define the $\SEP$ condition for a set of cones:
\Definition{Let $\Delta$ be a set of cones in a vector space $\R^{n}$ ($\Delta$ not necessarily being a fan) and $V:=\{v_{1}, ..., v_{p}\}$ a set of generating rays of $\Delta$ (that is, every cone of $\Delta$ is positively generated by a subfamily of $V$). We say that $\Delta$ \textbf{satisfies the $\boldsymbol{\SEP}$ condition} if for every cone $\sigma=\R_{\geqslant0}v_{i_{1}}+...+\R_{\geqslant0}v_{i_{k}} \in \Delta$ and every $i \in \{1, ..., p\}$, there exists $j \in \{i_{1}, ..., i_{k}\}$ such that $\sigma':=\R_{\geqslant0}v_{i_{1}}+...+\widehat{\R_{\geqslant0}v_{j}}+...+\R_{\geqslant0}v_{i_{k}}+\R_{\geqslant0}v_{i}\in\Delta$. }

\Lemme{[Cocompactness]\label{equivCond1}Let $\Delta$ be a simplicial fan in $\R^{n}$ (i.e. all its cones are simplicial) and $E \cong \R^{k}$ be a linear subspace of $\R^{n}$. Suppose that the quotient map $\pi: \R^{n} \rightarrow \R^{n}/E$ is injective on $|\Delta|$. Consider the set $\Delta_{\max}$ of all cones of $\Delta$ of maximal dimension. Then $\pi(\Delta)$ is complete if and only if the cones of $\Delta_{\max}$ are of dimension $n-k$ and $\Delta_{\max}$ satisfies the $\SEP$ condition.}
\Preuve{The assertion is clear.}


\subsubsection{From $\LVMB$ data to toric data}~\\
In this subsection we prove part $i)$ of theorem~\ref{BosioIsToric}. We see how one recovers toric information from an $\LVMB$ datum. 

Let $(\mathcal{L}, \mathcal{E}_{m,n})$ be an $\LVMB$ datum and $U$ the corresponding open subset of $\Pn$ (see section~\ref{BosioConstr}). We want to find a pair $(E,\Delta)$ (where $E$ is a $2m$-dimensional subspace of $\R^{n}$ and $\Delta$ is a fan of $\R^{n}$) satisfying conditions $a)$ and $b)$ of theorem~\ref{BosioIsToric}. Consider the fan $\Delta_{\Pn}$ in $\R^{n}$ defining the toric manifold $\Pn$: its rays are generated by the vectors $\{e_{1}, ..., e_{n}, -(e_{1}+...+e_{n})\}$ where $(e_{1},...,e_{n})$ is the canonical basis of $\R^{n}$. We call $e_{0}:=-(e_{1}+...+e_{n})$ and to each $P \in \mathcal{E}_{m,n}$ we associate the simplicial cone $\sigma_{P}$ of dimension $n-2m$ generated by the $n-2m$ vectors $e_{i}$ satisfying $i \not\in P$. A direct computation shows that the open set $U$ is the toric manifold given by the subfan of $\Delta_{\Pn}$ consisting of the cones $\{\sigma_{P}, P\in\mathcal{E}\}$ and their faces. Call this fan $\Delta$.

Now that we have found the fan $\Delta$, we must detect the subspace $E\cong \R^{2m}\subset \R^{n}$. For this, we need a preliminary lemma. Since $\C^{m}$ acts on $\Pn$ by $Z.[z_{0}: ... : z_{n}]= [\exp(\ell_{0}(Z))z_{0}:...:\exp(\ell_{n}(Z))z_{n}]$, we see $\C^{m}$ as a closed subgroup of $(\C^{*})^{n}$ by $$Z \mapsto (\exp(\ell_{1}(Z)-\ell_{0}(Z)), ..., \exp(\ell_{n}(Z)-\ell_{0}(Z))).$$

\Lemme{The intersection of the two subgroups $\C^{m}$ and $(\S^{1})^{n}$ of $(\C^{*})^{n}$ is trivial.}
\Preuve{
Choose $P\in \mathcal{E}_{m,n}$. An element $Z$ of $\C^{m} \cap (\S^{1})^{n}$ must satisfy
\[\ell_{k}(Z)-\ell_{0}(Z) \in i\R, \text{ for all $k \in P$;}\]
therefore we have $\Re(\ell_{0}(Z))=\Re(\ell_{k}(Z))$ for all $k \in P$. This means that $Z=0$ since $\{\ell_{k}, k \in P\}$ is a $\R$-affine basis for $(\C^{m})^{*}$.
}

Now we define $E\cong \R^{2m}\subset \R^{n}$ to be the image of $\C^{m}\subset (\C^{*})^{n}$ by the $\ord$ map:
$$\begin{array}{rrcl}\ord : & (\C^{*})^{n} & \rightarrow & \R^{n}\\
 & (z_{1}, ..., z_{n}) & \mapsto & (-\log|z_{1}|, ..., -\log|z_{n}|),\end{array}$$ which is injective on $\C^{m}$ by the previous lemma. A quick computation shows that a basis of this linear subspace of $\R^{n}$ is given by the following $2m$ vectors:

\begin{equation}\label{baseR2m}
X_{k}=\left(\begin{array}{c}\Re(\ell_{1,k}-\ell_{0,k}) \\ \vdots \\ \Re(\ell_{n,k}-\ell_{0,k})\\\end{array}\right)\!, Y_{k}=\left(\begin{array}{c}\Im(\ell_{1,k}-\ell_{0,k}) \\ \vdots \\ \Im(\ell_{n,k}-\ell_{0,k})\\\end{array}\right)\!, k\in\{1, ..., m\},
\end{equation} 
where $\ell_{i,k}$ is the $k$-th component of $\ell_{i}$.

We now have a pair $(E,\Delta)$ so we have to prove that they satisfy the needed properties. For this, first notice that we have the following commutative diagram: 
\begin{equation}
\label{diagrammeCommutatif}
\xymatrix{
U=X_{\Delta} \ar[rrrr]^{\displaystyle (\S^{1})^{n}}\ar[dd]_{\displaystyle  \C^{m}} & & & & \operatorname{\mathcal{M}c}(\Delta) \ar[dd]^{\displaystyle E= \ord(\C^{m})\cong\R^{2m}}\\
& & & & \\
X 
 \ar[rrrr]_{\displaystyle q((\S^{1})^{n}) \cong(\S^{1})^{n}}& & & & \operatorname{\mathcal{M}c}(\pi(\Delta)),
}\end{equation}
where $\ord : (\C^{*})^{n} \rightarrow (\C^{*})^{n}/(\S^{1})^{n}$ and $q : (\C^{*})^{n} \rightarrow (\C^{*})^{n}/\C^{m}$ are the quotient maps.
Let us call $\pi : \R^{n} \rightarrow \R^{n}/E\cong \R^{n-2m}$ the projection map with respect to $E$.
The fact that the quotient $\Mc(\Delta)/E$ is homeomorphic to $\Mc(\pi(\Delta))$ is easily checked. Indeed, the homeomorphism is the map $\widetilde{\pi}$ sending the orbit of $y+\infty\cdot\sigma$ to $\pi(y) +\infty\cdot\pi(\sigma)$. Bijectivity is a consequence of the injectivity of $\pi$ on $|\Delta|$. To prove it is open, consider an open set $U$ of $\Mc(\Delta)/E$ and a point $\widetilde{\pi}(x)=\widetilde{x}\in\widetilde{\pi}(U)$. Then, there are $\varepsilon>0$ and $w$ such that $U_{\varepsilon,w}(x)/E$ (see definition~\ref{topoMc}) is a neighbourhood of $x$ contained in $U$, and for $\delta>0$ small enough, $U_{\delta,\pi(w)}(\pi(x))$ is a subset of $\widetilde{\pi}( U_{\varepsilon,w}(x)/E )$. To prove continuity, we use a similar reasoning.~\\

Since the group $(\S^{1})^{n}$ is compact, the properness of the action of $\C^{m}$ on $U$ is equivalent to the properness of the action of $E$ on $\operatorname{\mathcal{M}c}(\Delta)$, which in turn is equivalent to the fact that $\pi$ is injective on $|\Delta|$ according to lemma~\ref{equivCond2}. On the other hand, $\operatorname{\mathcal{M}c}(\pi(\Delta))$ is compact because $X$ is, and the completeness of the fan $\pi(\Delta)$ is now a consequence of proposition~\ref{compactComplet}. By lemma~\ref{equivCond1}, one now can see that conditions $a)$ and $b)$ of theorem~\ref{BosioIsToric} are satisfied by the pair $(E,\Delta)$, hence part $i)$ of this theorem is proved.


~\\
\subsubsection{\label{ToricToBosio}From toric data to $\LVMB$ data.}~\\
We now prove part $ii)$ of theorem~\ref{BosioIsToric}. Suppose we are given a subfan $\Delta$ of $\Delta_{\Pn}$ and a linear subspace $E\cong \R^{2m} \subset \R^{n}$ satisfying both conditions $a)$ and $b)$ of theorem~\ref{BosioIsToric}. In order to recover an $\LVMB$ datum, we first choose $U = X_{\Delta}$. For $\sigma=\R_{\geqslant0}e_{i_{1}}+ ...+\R_{\geqslant0}e_{i_{n-2m}} \in \Delta$ a cone of maximal dimension $n-2m$, define $P_{\sigma}:=\{0,...,n\}\setminus\{i_{1}, ..., i_{n-2m}\}$ and $\mathcal{E}_{m,n} := \{P_{\sigma}, \sigma\in\Delta\text{ of dimension }n-2m\}$. It is clear that $U$ is the open set corresponding to $\mathcal{E}_{m,n}$ (as defined in section~\ref{BosioConstr}). 

We now have to see how one recovers the set of linear forms $\mathcal{L}$ and check if the pair $(\mathcal{L}, \mathcal{E}_{m,n})$ is an $\LVMB$ datum.

To do this, we first choose a basis for $E\cong \R^{2m}\subset \R^{n}$ which we write as a matrix $A = (a_{i,j}) \in \mathfrak{M}_{n,2m}(\R)$. Then we define $n+1$ vectors of $\C^{m}$ by taking each of the $n$ rows of $A$ and sending them to $\C^{m}$ via the map $(x_{1}, ..., x_{2m})\mapsto (x_{1}, ..., x_{m})+i(x_{m+1}, ..., x_{2m})$ (call these vectors $\ell_{1}, ..., \ell_{n}$), along with the zero vector of $\C^{m}$. Call this vector $\ell_{0}$ and set $\mathcal{L}=\{\ell_{0}, ..., \ell_{n}\}$. 

We now consider the pair $(\mathcal{L},\mathcal{E}_{m,n})$. According to lemma~\ref{equivCond2} we know that the action of $\C^{m}$ is proper on $X_{\Delta}$, meaning that the imbrication condition is satisfied. Since $\pi(\Delta)$ is complete, lemma~\ref{equivCond1} tells us that the set of all $n-2m$-dimensional cones of $\Delta$ satisfies the $\SEP$, hence $\mathcal{E}_{m,n}$ also does. Finally, the two Bosio conditions are satisfied, i.e. the pair $(\mathcal{L},\mathcal{E}_{m,n})$ is an $\LVMB$ datum.\\

\subsection{A correspondence between Bosio and toric data.}~\\
Let $(\mathcal{L},\mathcal{E}_{m,n})$ be an $\LVMB$ datum. For the proof of theorem~\ref{BosioIsToric}, we associated to this datum a unique pair $(E,\Delta)$ where $E$ is a $2m$-dimensional linear subspace of $\R^{n}$ and $\Delta$ is a fan of $\R^{n}$. 

As one can see from the proof of theorem~\ref{BosioIsToric}, it is possible that two different $\LVMB$ data give the same pair $(E,\Delta)$. We discuss this fact now.

Let $(\mathcal{L},\mathcal{E}_{m,n})$ and $(\mathcal{L'},\mathcal{E'}_{m',n'})$ be two $\LVMB$ data and let $(E, \Delta)$, $(E', \Delta')$ be the two associated pairs given by theorem~\ref{BosioIsToric}. If $(E, \Delta)=(E', \Delta')$, one has $n=n'$ and $m=m'$ because the dimensions of $E$ and $E'$ and their ambient spaces are equal respectively. Also, one sees that $\mathcal{E}_{m,n}=\mathcal{E'}_{m',n'}$ because it is clear from the proof of theorem~\ref{BosioIsToric} that for a subfan $\Delta$ of the fan of $\Pn$ there is a unique corresponding set $\mathcal{E}_{m,n}$. On the other hand, the fact that $E=E'$ does not imply $\mathcal{L}=\mathcal{L}'$, but it tells us that there exists a real affine automorphism of $\R^{2m}\cong(\C^{m})^{*}$ sending each $\ell_{i} \in\mathcal{L}$ to an element $\ell'_{i'} \in\mathcal{L}'$. 

Call $\mathcal{P}$ the set of pairs $(E,\Delta)$ such that $E$ and $\Delta$ satisfy conditions $a)$ and $b)$ of theorem~\ref{BosioIsToric}. Define an equivalence relation $\approx$ on the set of $\LVMB$ data, setting $(\mathcal{L},\mathcal{E}_{m,n}) \approx (\mathcal{L'},\mathcal{E'}_{m',n'})$ if and only if $\mathcal{E}_{m,n}=\mathcal{E'}_{m',n'}$ and there exists a real affine automorphism of $\R^{2m}\cong(\C^{m})^{*}$ such that its restriction to $\mathcal{L}$ is a bijection with $\mathcal{L}'$. 
The discussion above leads to the following statement:
\Proposition{There is a bijective correspondence
$$\{\LVMB \text{ data}\}/\approx~~~~ \longleftrightarrow \mathcal{P}.$$}

A natural question is now to ask when two $\LVMB$ data give the same manifolds (up to biholomorphism). We discuss this in section~\ref{biholLVMB} below.

\section{Detecting $\LVM$ data among $\LVMB$ data}
In \cite{Meersseman:2000aa}, Meersseman gave a method of construction of compact complex manifolds, called $\boldsymbol{\LVM}$ \textbf{manifolds}. Bosio shows that these manifolds can be obtained by his construction and gives a criterion to detect when an $\LVMB$ datum leads to an $\LVM$ manifold. This is proposition~1.3 in \cite{Bosio:2001aa} which we recall now. 

Let $(\mathcal{L},\mathcal{E}_{m,n})$ be an $\LVMB$ datum and define $\mathcal{O}$ to be the set of points in $(\C^{m})^{*}$ (the dual space of $\C^{m}$) which are not in the convex hull of any family of $2m$ elements of $\mathcal{L}$. 

\Proposition{\label{critereLVMBosio}Given an $\LVMB$ datum $(\mathcal{L},\mathcal{E}_{m,n})$, one obtains an $\LVM$ manifold if and only if there is a bounded connected component $O$\! of $\mathcal{O}$ such that $\mathcal{E}_{m,n}$ is the collection of subsets $P$ of $\{1,...,n\}$ having $2m+1$ elements, with the property that the convex envelope of $\mathcal{L}_{P}$ contains $O$.}
\Definition{We say that an $\LVMB$ datum is an $\boldsymbol{\LVM}$ \textbf{datum} if it satisfies the previous condition.}

Notice that proposition~\ref{critereLVMBosio} can be rephrased as follows:
\Proposition{An $\LVMB$ datum $(\mathcal{L},\mathcal{E}_{m,n})$ is an $\LVM$ datum if and only if $$\bigcap_{P\in\mathcal{E}_{m,n}} \mathring{C}_{P} \neq \emptyset,$$
where $ \mathring{C}_{P} $ is the interior of the convex envelope of $\mathcal{L}_{P}$.}

In \cite{Cupit-Foutou:2007aa}, Cupit-Foutou and Zaffran give another characterization of $\LVM$ data among $\LVMB$ data with a supplementary assumption, called ``condition (K)'', see proposition~3.2 in \cite{Cupit-Foutou:2007aa}. (We say that an $\LVMB$ datum $(\mathcal{L},\mathcal{E}_{m,n})$ satisfies condition (K) if there exists an affine automorphism of $(\C^{m})^{*}$ which maps each $\ell_{i}\in\mathcal{L}$ to a vector with integer coefficients.) Theorem~\ref{lvmbEquivPolytopal} below can be seen as a generalization of their result. Before we can state it and give a proof, we need to recall some preliminary definitions and statements:

\Definition{A fan $\Delta$ is called \textbf{strongly polytopal} (for short, in the following, \textbf{polytopal}) if there exists a polytope $P$ containing $0$ in its interior, such that $\Delta$ is the set of cones generated by the faces of $P$.}

\Exemple{Every complete fan in $\R^{2}$ is polytopal.}

Shephard gives in \cite{Shephard:1971aa} a criterion for a fan to be polytopal, we recall it now. 

\Definition{(See for instance \cite{Grunbaum:2003aa} or \cite{Ewald:1996aa}.) Let $X=(x_{1}, ..., x_{r}) \in \R^{n}$ be a family of vectors and $A$ the matrix whose columns are the elements of $X$ (hence, $A$ is a matrix with $r$ columns and $n$ rows). Also assume that $\dim \Aff X = r$ where $\Aff X$ is the affine hull of $X$. Let $(\alpha_{1}, ..., \alpha_{r-n})$ be a basis of the kernel of $A$ and $B$ the $(r-n) \times n$ matrix whose columns are the $\alpha_{i}'s$. The family $\overline{X}=(\overline{x}_{1}, ..., \overline{x}_{r}) \subset \R^{r-n}$ of vectors given by the rows of $B$ is called \textbf{a linear transform} of $X$. }
Notice that we can not talk about \emph{the} linear transform of a family $X$, since the construction depends on the choice of the basis of the kernel of $A$. 

We have the following lemma (see for instance \cite{Ewald:1996aa}):
\Lemme{Let $X=(x_{1}, ..., x_{r})$ be a family of vectors as in the previous definition and $\overline{X}=(\overline{x}_{1}, ..., \overline{x}_{r})$ a linear transform of $X$. Notice that $X$ is a linear transform of $\overline{X}$. Then $\textstyle \sum\limits_{\vphantom{I}i=1}^{r} \overline{x}_{i}=0$ (resp. $\textstyle  \sum\limits_{\vphantom{I}i=1}^{r} x_{i}=0$) if and only if the vectors $x_{i}$ (resp. $\overline{x}_{i}$) all belong to a hyperplane $H \subset \R^{n}$ (resp. $\subset \R^{r-n}$) which does not contain $0$. }

\Definition{Let $X=(x_{1}, ..., x_{r})$ be a family of vectors in $\R^{n}$ which positively span $\R^{n}$ (i.e. $\R_{\geqslant0}x_{1} + ... + \R_{\geqslant0}x_{r} = \R^{n}$). Choose a suitable family $\lambda_{i}>0$ such that $\sum \lambda_{i}x_{i} = 0$ (this is always possible, see \cite{Shephard:1971aa}, p.~258). Choose a linear transform of the family $(\lambda_{1} x_{1}, ..., \lambda_{r}x_{r})$ such that the last coordinate is always equal to $1$ and call such a family a \textbf{Shephard transform} of $X$. Denote it by $\widehat{X}=(\widehat{x}_{1}, ..., \widehat{x}_{r})$. }

Now consider a complete fan $\Delta$ in $\R^{n}$ and assume that $X=(x_{1}, ..., x_{r})$ is a family of vectors which generate the rays of $\Delta$. If $\sigma = \R_{\geqslant 0}x_{i_{1}}+\cdots +\R_{\geqslant 0} x_{i_{p}}$ is a cone of $\Delta$ of (maximal) dimension $n$, we denote by $\mathring{C}(\sigma)$ the relative interior of the convex envelope of $\widehat{X}\setminus\{\widehat{x}_{i_{1}}, ..., \widehat{x}_{i_{p}}\}$. Then we have the following:

\Theoreme{[Shephard's criterion]\label{critShephard} With the notations above, the fan $\Delta$ is polytopal if and only if $$\bigcap_{\sigma \in \Delta_{\max}} \mathring{C}(\sigma)\neq\emptyset,$$ where $\Delta_{\max}$ is the set of all cones of $\Delta$ of maximal dimension $n$.}

For the previous definition and theorem we refer to \cite{Shephard:1971aa}. We are now able to prove the following:

\Theoreme{\label{lvmbEquivPolytopal}Let $(\mathcal{L},\mathcal{E}_{m,n})$ be an $\LVMB$ datum and $(E,\Delta)$ its associated pair given by theorem~\ref{BosioIsToric}. Then, $(\mathcal{L},\mathcal{E}_{m,n})$ is an $\LVM$ datum if and only if the projection by $E$ of the fan $\Delta$ is polytopal.}

\Preuve{Form a basis of $\R^{n}$ with first $2m$ vectors being the vectors $X_{k}, Y_{k}$ defined in equation~(\ref{baseR2m}) which generate $E$ and, for last $n-2m$ vectors, any of the canonical basis. 
Call $\mathcal{B}'$ this basis, and $\mathcal{B}$ the canonical basis of $\R^{n}$:
$$\mathcal{B}'=(X_{1}, ..., Y_{2m}, e_{i_{1}}, ..., e_{i_{n-2m}}).$$ 
We have the direct sum decomposition $\R^{n} = E\oplus F$ where $F\cong \R^{n-2m}$. The projection onto $F$ with respect to $E$ is given by the matrix $\Pi:= J_{n,m} P^{-1}$ where $P$ is the invertible matrix whose columns are the vectors of $\mathcal{B}'$ and 
$$J_{n,m} = 
\left(\begin{array}{c|c}0_{n-2m,2m} & I_{n-2m} \end{array}\right).$$

(Here, $0_{n-2m,2m}$ is the zero matrix with $n-2m$ rows and $2m$ columns, and $I_{n-2m}$ is the identity matrix of size $n-2m$.) We now see that the vectors $e_{1}, ..., e_{n}$ and $e_{0}=-(e_{1}+...+e_{n})$ are sent by $\Pi$ respectively to the $n$ columns of $\Pi$ and the opposite of their sum. 
~\\

We now separate cases, depending on the number of indispensable $\ell_{i}$'s (note that there can be between $0$ and $2m$ indispensable $\ell_{i}$'s).~\\

First, assume that there is no indispensable $\ell_{i}$. It is readily seen that the vectors $(\Pi.e_{0}, \Pi.e_{1}, ..., \Pi.e_{n})$ generate the rays of the (complete) projected fan $\pi(\Delta)$ of $F\cong \R^{n-2m}$. 
It is also straightforward to check that a Shephard transform of this family of vectors is given by the vectors $\widehat{\ell}_{i} := (\ell_{i}, 1) \in \R^{2m+1}$, for $i=0, ..., n$. Now, the use of both Shephard and Bosio criteria (theorem~\ref{critShephard} and proposition~\ref{critereLVMBosio}) gives the equivalence between the polytopality of $\pi(\Delta)$ and the fact that $(\mathcal{L},\mathcal{E}_{m,n})$ is $\LVM$.~\\

Assume now that there is one indispensable $\ell_{i}$, say $\ell_{0}$. First, notice that it implies that the $1$-dimensional cones of the fan $\Delta$ are generated by $e_{1}, ..., e_{n}$ (and not $e_{0}$). Then, in order to caracterize the polytopality of $\pi(\Delta)$, one must compute a Shephard transform of $(\Pi.e_{1}, ..., \Pi.e_{n})$, which we do in two steps. First, we project the vectors $\widehat{\ell}_{i}$ (for $i=1,..,n$) with respect to $\R.\widehat{\ell}_{0}$. This procedure leads to a linear transform of $(\Pi.e_{1}, ..., \Pi.e_{n})$. Recall that the vectors $\widehat{\ell}_{i}$ belong to the affine hyperplane $$H:= \{(x_{1}, ..., x_{2m+1})\in\R^{2m+1}~|~x_{2m+1}=1\},$$
 so we can define $H_{0}$ to be the hyperplane $H-\widehat{\ell}_{0}$, then one sees that $H_{0}$ is a linear subspace of $\R^{2m+1}$, and that $\widehat{\ell}_{0} \not \in H_{0}$. Call $\pi_{0} : \R^{2m+1} \rightarrow H_{0}$ the projection with respect to $\R.\widehat{\ell}_{0}$. Endow $H$ with a vector space structure by choosing $\widehat{\ell}_{0}$ as its origin, then the restriction of $\pi_{0}$ between $H$ and $H_{0}$ is an isomorphism. 
Remark that if $\ell_{0}$ is indispensable, it can not be an element of the convex envelope of $\{\ell_{1}, ..., \ell_{n}\}$. Indeed, since $(\Pi.e_{1}, ..., \Pi.e_{n})$ positively spans $\R^{n-2m}$, one has $\pi_{0}(\widehat{\ell}_{0}) \not \in \operatorname{conv}(\pi_{0}(\widehat{\ell}_{1}), ..., \pi_{0}(\widehat{\ell}_{n}))$ (see \cite{Shephard:1971aa}). Hence, there exists a hyperplane $H'_{0}$ of $H_{0}$ containing $2m$ elements of $\{\pi_{0}(\widehat{\ell}_{1}), ..., \pi_{0}(\widehat{\ell}_{n})\}$, separating $\pi_{0}(\widehat{\ell}_{0})=0\in H_{0}$ and $\{\pi_{0}(\widehat{\ell}_{1}), ..., \pi_{0}(\widehat{\ell}_{n})\}$.  See figure~\ref{fig2} for a picture of the situation. 

\begin{figure}[h]
\includegraphics[scale=1]{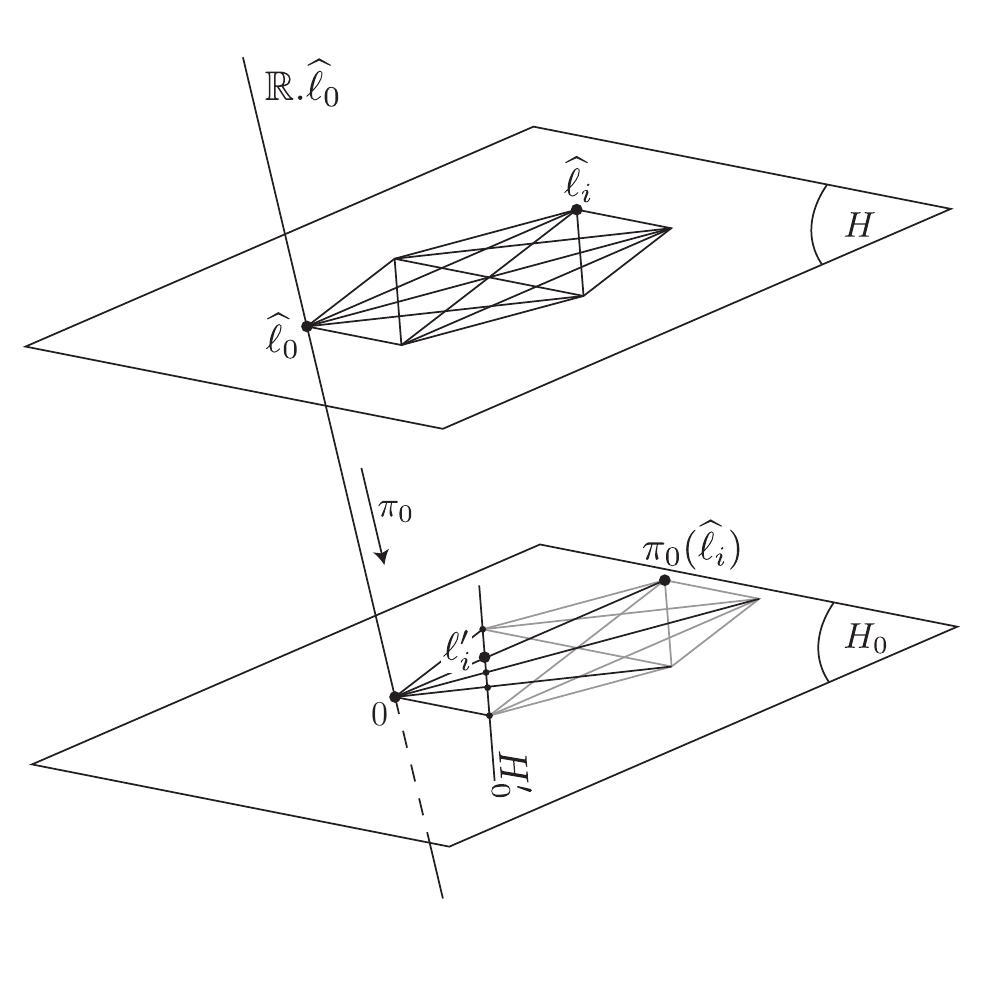}
\caption{\label{fig2}The projection $\pi_{0}$}
\end{figure}

The second step we have to do in order to obtain a Shephard transform of $(\Pi.e_{1}, ..., \Pi.e_{n})$ is to multiply each vector of $(\pi_{0}(\widehat{\ell}_{1}), ..., \pi_{0}(\widehat{\ell}_{n}))$ by suitable scalars such that they all belong to $H'_{0}$. Call $\mathcal{L}':=(\ell'_{1}, ..., \ell'_{n})$ the vectors obtained that way. Remark that any family of $2m$ vectors of this set forms an affine basis for $H'_{0}$. For $P = (0, i_{1}, ..., i_{2m}) \in \mathcal{E}_{m,n}$, call $P':=(i_{1}, ..., i_{2m})$ and denote by $\mathcal{E}'_{m,n}$ the set $\{P'~|~P \in\mathcal{E}_{m,n}\}$. It is clear that $\mathcal{E}'_{m,n}$ still satisfies Bosio's $\SEP$ (with $\{1, ..., n\}$ instead of $\{0,...,n\}$). We show that to a point in the convex envelope of a family of points $\mathcal{L}_{P}$ it is possible to associate a point in the convex envelope of $\mathcal{L}'_{P'}$ and reciprocally; this fact proves that the imbrication condition is also satisfied by $(\mathcal{L}',\mathcal{E}'_{m,n})$. Assume $x_{0} \in \operatorname{conv}(\mathcal{L}'_{P'})$ for $P' \in \mathcal{E'}_{m,n}$, then the point $\frac{1}{2}x_{0}$ is an element of $\operatorname{conv}(\mathcal{L}_{P})$. Reciprocally if $x_{0} \in \operatorname{conv}(\mathcal{L}_{P})$ (for $P\in \mathcal{E}_{m,n}$), there exists a unique $\lambda > 0$ (which does not depend on $P$) such that $\lambda x_{0} \in H'_{0}$. Write $x_{0} = \displaystyle\sum_{j=1}^{2m}\lambda_{j} \beta_{i_{j}}\ell'_{i_{j}}$ (where $\beta_{i_{j}}\ell'_{i_{j}}=\pi_{0}(\widehat{\ell_{i_{j}}})$), and $\lambda x_{0} = \displaystyle \sum_{j=1}^{2m}\alpha_{j}\ell'_{i_{j}}$ with $\sum_{j}\alpha_{j}=1$. Since any subfamily of $\mathcal{L}'$ of $2m$ elements is an $\R$-affine basis of $H'_{0}$, one has $\lambda \lambda_{j} \beta_{i_{j}}= \alpha_{j}$, hence $\alpha_{j}\geqslant 0$ for all $j$, that is $\lambda x_{0} \in \operatorname{conv}(\mathcal{L}_{P'}')$.

An important consequence of this remark is also that $\displaystyle \bigcap_{P'\in\mathcal{E}'_{m,n}} \!\!\!\! \operatorname{conv}(\mathcal{L}'_{P'})\neq \emptyset$ if and only if $\displaystyle \bigcap_{P\in\mathcal{E}_{m,n}}\!\!\!\! \operatorname{conv}(\mathcal{L}_{P}) \neq \emptyset$.

Now, the consecutive use of Shephard's criterion and Bosio's criterion applied to the Shephard transform that we have computed gives us that $\pi(\Delta)$ is polytopal if and only if $(\mathcal{L},\mathcal{E}_{m,n})$ is $\LVM$.

The same reasoning now applies if there are up to $2m-1$ indispensable $\ell_{i}$'s, say $\ell_{0}, ..., \ell_{k}$ by applying the previous method ``recursively'': to compute a Shephard transform of $(\Pi.e_{k+1}, ..., \Pi.e_{n})$, one can compute a Shephard transform of $(\Pi.e_{1}, ..., \Pi.e_{n})$ just as above, then project it with respect to $\R.\widehat{\Pi.e_{1}}$ and multiply by the proper scalars, which gives a Shephard transform of $(\Pi.e_{2}, ..., \Pi.e_{n})$ and so on. Call $(\mathcal{L}^{(j)}, \mathcal{E}^{(j)}_{m,n})$ the families obtained at each step (for $j=1,...,k+1$) and $(\mathcal{L}^{(0)}, \mathcal{E}^{(0)}_{m,n}):=(\mathcal{L}, \mathcal{E}_{m,n})$. Notice (just as above) that for $j=1,...,k+1$, the intersection $\displaystyle \bigcap_{P\in\mathcal{E}^{(j)}_{m,n}}\!\!\!\! \operatorname{conv}(\mathcal{L}^{(j)}_{P})$ is non-empty if and only if $\displaystyle \bigcap_{P\in\mathcal{E}^{(j-1)}_{m,n}}\!\!\!\! \operatorname{conv}(\mathcal{L}^{(j-1)}_{P})$ is non-empty. Then, using Shephard and Bosio criteria again, we obtain the result.

~\\Finally, assume there are $2m$ indispensable items, 
say $\{0, ..., 2m-1\}$ for simplicity. The $\SEP$ condition then forces that $$\mathcal{E}_{m,n} = \{(0, ..., 2m-1, j) ~|~j \in \{2m, ..., n\}\}.$$ The second Bosio condition (imbrication) implies that the hyperplane $H$ containing $\{\ell_{0}, ..., \ell_{2m-1}\}$ is a supporting hyperplane for the polytope obtained as the convex envelope of $\mathcal{L}$. As it can be easily seen, there is exactly one bounded connected component $O$ of $\mathcal{O}$ such that $\Aff(\overline{O}\cap H)=H$, and this component lies in the intersection of all the convex envelopes of $\{\mathcal{L}_{P}, P \in\mathcal{E}_{m,n}\}$, hence the $\LVMB$ datum we consider is $\LVM$ by proposition~\ref{critereLVMBosio}. On the other hand, the fact that there are $2m$ indispensable elements says that the projected fan has exactly $n+1-2m$ generating rays and each cone of maximal dimension is generated by $n-2m$ vectors. Such a complete fan in $\R^{n-2m}$ is always polytopal. }

\section{\label{biholLVMB}When are two $\LVMB$ manifolds biholomorphic?}
It is possible that two different $\LVMB$ data give rise to two biholomorphic manifolds, this is why we make a distinction between $\LVMB$ \emph{data} and \emph{manifolds}. 
In this section we prove that if an $\LVMB$ manifold $X_{n,m}$ built with an $\LVMB$ datum $(\mathcal{L},\mathcal{E}_{m,n})$ is boholomorphic to an $\LVM$ manifold, then the datum $(\mathcal{L},\mathcal{E}_{m,n})$ itself is $\LVM$. This problem was raised and partly answered by Cupit-Foutou and Zaffran in \cite{Cupit-Foutou:2007aa}. For the proof, we will use the criterion given by theorem~\ref{lvmbEquivPolytopal}.\\

We now study the case when two $\LVMB$ data $(\mathcal{L}_{1}, \mathcal{E}_{m_{1},n_{1}})$ and $(\mathcal{L}_{2}, \mathcal{E}'_{m_{2},n_{2}})$ give the same manifold (up to biholomorphism).\\

Let $X$ be an $\LVMB$ manifold. Let $A:=Aut_{\mathcal{O}}(X)^{o}$ be the identity component of the group of holomorphic automorphisms of $X$ which is a complex Lie group because $X$ is compact (this is a result of Bochner and Montgomery, see \cite{Bochner:1947aa}). \\

We first prove the following
\Lemme{Let $(\mathcal{L}_{1}, \mathcal{E}_{m_{1},n_{1}})$ and $(\mathcal{L}_{2}, \mathcal{E}'_{m_{2},n_{2}})$ be two $\LVMB$ data such that the corresponding $\LVMB$ manifolds are biholomorphic. Then, $m_{1}=m_{2}$ and $n_{1}=n_{2}$.}
\Preuve{Assume that $X_{1}$ and $X_{2}$ are two isomorphic $\LVMB$ manifolds obtained as quotient of two open sets $U_{1}\subset \P^{n_{1}}(\C)$ and $U_{2}\subset \P^{n_{2}}(\C)$ by two closed subgroups $H_{1}\cong\C^{m_{1}}\subset(\C^{*})^{n_{1}}$ and $H_{2}\cong \C^{m_{2}}\subset(\C^{*})^{n_{2}}$. Call $\varphi$ a biholomorphism between $X_{1}$ and $X_{2}$. It induces an isomorphism $$\varphi^{*} : Aut_{\mathcal{O}}(X_{2})^{o}=:A_{2} \rightarrow A_{1}:=Aut_{\mathcal{O}}(X_{1})^{o},$$ given by $\varphi^{*}(g) = \varphi^{-1}\circ g\circ \varphi$.

Then, $G_{i}:=(\C^{*})^{n_{i}}/H_{i}$ is a subgroup of $A_{i}$ and has an open dense orbit $G_{i}.x_{i} \subset X_{i}$ for $i=1,2$ (see \cite{Bosio:2001aa}, proposition~2.4). Call $\pi_{i} : (\C^{*})^{n_{i}} \rightarrow (\C^{*})^{n_{i}}/H_{i}$ the two quotient maps.
Now the subgroup $T_{1}:=\pi_{1}((\S^{1})^{n_{1}}) \subset G_{1}$ (resp. $T_{2}=\pi_{2}((\S^{1})^{n_{2}}) \subset G_{2}$) is contained in a maximal torus of $A_{1}$ (resp. $A_{2}$), 
say $\widetilde{T}_{1}$ (resp. $\widetilde{T}_{2}$). The two maximal tori $\widetilde{T}_{1}$ and $\varphi^{*}(\widetilde{T}_{2})$ of $A_{1}$ are conjugated (see, for instance, \cite{Brocker:1985aa}) so, up to conjugacy, we can assume that $T_{1}$ and $\varphi^{*}(T_{2})$ are contained in the same maximal torus $\widetilde{T}$ of $A_{1}$. 
Now we denote by $\widetilde{T}^{\C}$ the complexification of $\widetilde{T}$, and both $G_{1}$ and $\varphi^{*}(G_{2})$ are Lie subgroups of $\widetilde{T}^{\C}$ (because $T_{i}^{\C}=G_{i}$ for $i=1,2$).

Suppose we have $\dim_{\C}\widetilde{T}^{\C}> \dim_{\C} G_{1}=\dim_{\C} X_{1}$, then the action of $A_{1}$ is not effective anymore (in this case indeed, there exists a non-trivial element $f\in\widetilde{T}^{\C}$ such that $f(x_{1})=x_{1}$ and since $\widetilde{T}^{\C}$ is abelian (and contains $G_{1}$), $f$ is the identity map on all $G_{1}.x_{1}$, hence on all $X_{1}$ by analytic  continuation), a contradiction. Hence we necessarily have $\dim_{\C}\widetilde{T}^{\C}= \dim_{\C} G_{1}=\dim_{\C} \varphi^{*}(G_{2})$ and $G_{1} = \varphi^{*}(G_{2})$. In particular, these two Lie groups have isomorphic fundamental groups, $\Z^{n_{1}}$ and $\Z^{n_{2}}$ respectively, which leads to $n_{1}=n_{2}$ and since $n_{1}-m_{2}=n_{2}-m_{2}$, we also have $m_{1}=m_{2}$.}

Now we prove the following 
\Theoreme{\label{questionByZaffran}Let $(\mathcal{L}_{1}, \mathcal{E}_{m,n})$ and $(\mathcal{L}_{2}, \mathcal{E}'_{m',n'})$ be two $\LVMB$ data giving two biholomorphic $\LVMB$ manifolds, then $n=n'$, $m=m'$ and $(\mathcal{L}_{1}, \mathcal{E}_{m,n})$ is an $\LVM$ datum if and only if $(\mathcal{L}_{2}, \mathcal{E}'_{m,n})$ is an $\LVM$ datum.}
\Preuve{We keep the notations of the previous proof. Up to conjugacy, we can assume that $(\varphi^{-1})^{*}(G_{1})=G_{2}$ and $\varphi(G_{1}.x_{1})=G_{2}.x_{2}$, where $G_{i}.x_{i}$ is the open dense orbit of the action of $G_{i}$ on $X_{i}$. Call $(E,\Delta)$ and $(E',\Delta')$ the toric data associated respectively to these two $\LVMB$ data and let $\pi : \R^{n} \rightarrow \R^{n}/E$ and $\pi' : \R^{n} \rightarrow \R^{n}/E'$ be the corresponding projections.

The maximal compact tori of $G_{1}$ and $G_{2}$ are isomorphic to $(\S^{1})^{n}$, hence the restriction of $(\varphi^{-1})^{*}$ to these maximal tori is given by a matrix $A \in \operatorname{GL}(n,\Z)$, and it extends to an isomorphism of $(\C^{*})^{n}$ which we call $\widetilde{\varphi}$.

As one can check by passing to the Lie algebras, the following diagram is commutative: 
 \begin{equation}\label{diagram}\xymatrix{
    (\C^{*})^{n} \ar[r]^{\displaystyle \widetilde{\varphi}} \ar@{^{(}->}[d]_{\displaystyle\iota_{1}}& (\C^{*})^{n}\ar@{^{(}->}[d]^{\displaystyle\iota_{2}}\\
    U_{1}  \ar[d]_{\displaystyle p_{1}} & U_{2} \ar[d]^{\displaystyle p_{2}} \\
    X_{1} \ar[r]_{\displaystyle \varphi} & X_{2},
  }\end{equation}
where $\iota_{1}$ and $\iota_{2}$ are the canonical injections. 

By using the Riemann extension theorem, one sees that the application $p_{2} \circ\iota_{2}\circ \widetilde{\varphi}  = \varphi \circ p_{1}\circ\iota_{1}$ extends on $U_{1}$, that is, we can extend $\widetilde{\varphi}$ to a toric application between $U_{1}$ and $U_{2}$, so it induces a (bijective) linear map from $\R^{n}$ to itself (given by the matrix $A$), which maps the fan $\Delta$ of $\R^{n}$ defining the toric manifold $U_{1}$ to the fan $\Delta'$ of $\R^{n}$ defining $U_{2}$ (\cite{Oda:1988aa}, theorem~1.13). The commutativity of diagram~(\ref{diagram}) implies that this linear map induces a linear map from $\R^{n}/E\cong \R^{n-2m}$ to $\R^{n}/E'\cong \R^{n-2m}$ sending $\pi(\Delta)$ to $\pi'(\Delta')$. This means that these two projected fans are simultaneously polytopal and the conclusion now comes from theorem~\ref{lvmbEquivPolytopal}. }


\section{A generalization of Bosio's construction}
\subsection{Detecting all open subsets with a compact quotient}
The following proposition tells us that Bosio's construction is the most general in this context, in the sense that any open subset $U \subset \Pn$ with an action of $\C^{m}$ such that the quotient is a complex compact manifold is in fact given by a subfan of $\Delta_{\Pn}$.

\Proposition{\label{ouvertsToriques}Let $G\cong \C^{m} \subset (\C^{*})^{n}$ be a closed subgroup and $U\subset \Pn$ be an open subset such that $U/G$ is a compact manifold. Then $U$ is stable by $(\C^{*})^{n}$, therefore given by a subfan of the fan of $\Pn$.}
\Preuve{Let $\pi$ be the quotient map $\pi : U\rightarrow U/G$ and let $v_{1}, ..., v_{n}$ be a basis of the Lie algebra of $(\C^{*})^{n}$ such that $v_{1}, ..., v_{m}$ is a basis of the Lie algebra of $G$. We define the vector fields $v^{*}_{1}, ..., v^{*}_{n}$ on $U/G$ by $v_{i}^{*}(x):=d\pi_{y}(v_{i}(y))$ for $y\in\pi^{-1}(x)$. Since $U$ is stable by $G$, it is enough to show that $U$ is invariant by the action of $v_{m+1}, ... , v_{n}$. Assume this is not the case, say for instance that $U$ is not stable under the action of $v_{m+1}$: there exists a point $z_{0}\in U$, a holomorphic map $\gamma:\C\rightarrow \Pn$ and $t_{0}\in \C$ such that $\dot{\gamma}(0)=v_{m+1}(z_{0})$, $\gamma(0)=z_{0}$ and $\gamma(t_{0})\not\in U$. Call $\Omega:=\gamma^{-1}(U)$. Since $U/G$ is compact, the vector field $v^{*}_{m+1}$ is complete and $\pi(\gamma|_{\Omega})$ is an integral curve of $v^{*}_{m+1}$ passing through $\pi(z_{0})$ hence it can be extended to $\C$, a contradiction because then $\gamma(t_{0}) \in U$.}
%
%

\subsection{A generalization}
In light of the proof of proposition~\ref{ouvertsToriques}, where the property of $\Pn$ used is the fact that it is a compact toric variety, we can study a ``toric'' open subset of any compact toric manifold: 
\Theoreme{Let $\Delta$ be a finite rational fan in $\R^{n}$ and $E\cong \R^{2m}$ be a linear subspace of $\R^{n}$ such that: 
\begin{itemize}
\item[-] the projection map $\pi: \R^{n} \rightarrow \R^{n}/E\cong\R^{n-2m}$ is injective on $|\Delta|$,
\item[-] the fan $\pi(\Delta)$ is complete in $\R^{n}/E$, i.e. $|\pi(\Delta)| = \R^{n}/E$.
\end{itemize}

Define a closed subgroup $G \cong \C^{m}$ of $(\C^{*})^{n}$ in the same way as in section~\ref{ToricToBosio}. Then, the quotient $X_{\Delta} / G$ exists and it is a complex compact manifold.} 
\Preuve{Since $\Delta$ is finite we can construct a complete rational fan $\Delta'$ containing $\Delta$ as a subfan. The group $G$ acts properly and freely on $X_{\Delta}$ by lemma~\ref{equivCond2} hence we know that $X_{\Delta}/G$ is a complex manifold. Compactness is a consequence of the completeness of $\pi(\Delta)$.}

\section*{Acknowledgements}
I would like to thank Karl OELJEKLAUS for his constant support, helpful discussions, comments and careful reading of the paper.

\bibliographystyle{amsplain}
\bibliography{bibliographie}

\end{document}